\documentclass[a4paper, 10pt, conference]{ieeeconf}      

\IEEEoverridecommandlockouts                              
\usepackage{graphics} 
\usepackage{epsfig} 
\usepackage{mathptmx} 
\usepackage{times} 
\usepackage{amsmath} 
\usepackage{amssymb, verbatim}  
\usepackage{url}

\usepackage{color}
\usepackage{cite}

\newcommand{\I}{\textrm{I}}
\newcommand{\II}{\textrm{II}}
\newcommand{\III}{\textrm{III}}

\newtheorem{remark}{Remark}

\title{\LARGE \bf
A hierarchical time-splitting approach for \\ solving finite-time optimal control problems
}

\author{Georgios Stathopoulos$^{1}$, Tam{\'a}s Keviczky$^{2}$ and Yang Wang$^{3}$ 
\thanks{$^{1}$G. Stathopoulos is with Laboratoire d'Automatique, EPFL, CH-1015 Lausanne, Switzerland,
        {\tt\footnotesize georgios.stathopoulos@epfl.ch}}%
\thanks{$^{2}$T. Keviczky is with the Delft Center for Systems and Control, Delft University of Technology,
        Delft, CD 2628, The Netherlands,
        {\tt\footnotesize t.keviczky@tudelft.nl}}%
\thanks{$^{3}$Y. Wang is with Stanford University, Stanford, CA 94305, USA,
	{\tt\footnotesize yang1024@gmail.com}}%
}

\newcommand{\BEAS}{\begin{eqnarray*}}
\newcommand{\EEAS}{\end{eqnarray*}}
\newcommand{\BEQ}{\begin{equation}}
\newcommand{\EEQ}{\end{equation}}
\newcommand{\BIT}{\begin{itemize}}
\newcommand{\EIT}{\end{itemize}}

\newcommand{\reals}{{\mbox{\bf R}}}

\newcommand{\argmin}{\mathop{\rm argmin}}
\newcommand{\argmax}{\mathop{\rm argmax}}

\begin{document}

\maketitle
\thispagestyle{empty}
\pagestyle{empty}

\begin{abstract}


We present a hierarchical computation approach for solving finite-time optimal control problems using operator splitting methods.
The first split is performed over the time index and leads to as many subproblems as the length of the prediction horizon. Each
subproblem is solved in parallel and further split into three by separating the objective from the equality and inequality constraints
respectively, such that an analytic solution can be achieved for each subproblem. The proposed solution approach leads to a nested
decomposition scheme, which is highly parallelizable. We present a numerical comparison with standard state-of-the-art solvers,
and provide analytic solutions to several elements of the algorithm, which enhances its applicability in fast large-scale applications.

\end{abstract}

\section{INTRODUCTION}

Online optimization and optimal control methods are increasingly being considered for fast embedded applications,
where efficient, reliable, and predictable computations involved in calculating the optimal solutions are a necessity.
The potential use of optimal control in such embedded systems promises energy savings and more efficient resource usage,
increased safety, and improved fault detection. The range of application areas that can benefit from embedded optimization
include the mechatronics, automotive, process control and aerospace sectors \cite{EMBOCON}. The promise of unprecedented performance
and capabilities in these applications, which typically rely on large-volume, real-time embedded control systems,
has fueled recent research efforts towards fast and parallel optimization solvers. 

One of the main research directions aim at developing special-purpose optimization solvers that target typical control or estimation problems arising in
optimal control. Parallel solutions to systems of linear equations appearing in interior-point, and active set methods have been studied in
\cite{C2009,GBL2010,JCKL2011,WMN2011}. In this work we consider a quadratic finite-time optimal control problem for discrete-time systems with constrained linear dynamics, which appears in typical
model predictive control problems \cite{BBM2012}. We investigate and develop different parallelizable algorithms using operator splitting techniques \cite{PTVF2007,OSB2012}
that have recently shown great promise for speeding up calculations involved in computing optimal solutions with medium accuracy \cite{GS2012,KF2012,port_opt_bound}.
Our approach relies on a hierarchical splitting up of the specially structured finite-time optimal control problem. The first split is performed over the time
index and leads to as many subproblems as the length of the prediction horizon. Each subproblem can then be solved in parallel and further split into three
by separating the objective from the equality and inequality constraints respectively, such that an analytic solution can be achieved for each subproblem.
The proposed solution approach leads to a nested decomposition scheme, which is highly parallelizable. The proposed three-set splitting method does not only
solve the particular quadratic programs (QPs) that appear in the update steps of the time-splitting algorithm efficiently, but also provides a compact,
standalone alternative for solving generic QPs.

The paper is structured as follows. Section~\ref{sec:timesplitting} presents the main idea behind the time-splitting optimal control approach for parallel
computations, using the Alternating Direction Method of Multipliers and deriving the exact formulas required for each subproblem and update step. 
In Section~\ref{sec:threeset} we propose an alternative scheme for solving the QPs with general polyhedral constraints that arise in the time-splitting
update steps (or for any other generic QP). The two splitting schemes are combined in a hierarchical fashion in Section~\ref{sec:hierarchical}, and numerical
experiments are performed in Section~\ref{sec:results} to compare its performance with some advanced solvers in the literature. Section~\ref{sec:conclusions}
concludes the paper.

\section{TIME-SPLITTING OPTIMAL CONTROL} \label{sec:timesplitting}

\subsection{Problem Formulation}

We consider the following finite-time optimal control problem formulation that arises in typical model predictive control applications:
    \begin{subequations}   \label{eq:time_split_objective}
    \begin{eqnarray}
      \mbox{minimize} && \frac{1}{2}\sum_{t=0}^N \Big(x_t^T Q x_t + u_t^{T} R u_t\Big) \label{eq:time_split_objective_cost} \\
      \mbox{subject to} && (x_t,u_t) \in \mathcal{X}_t \times \mathcal{U}_t, \quad t = 0,\ldots,N  \label{eq:time_split_objective_constr} \\
	  \quad && x_{t+1} = A_t x_t+B_t u_t+c_t, t = 0,\ldots,N-1 \label{eq:time_split_objective_dyn}
    \end{eqnarray}
    \end{subequations}
where the decision variables are the states $x_t \in \reals^n$, and the inputs $u_t \in \reals^m$ of the system for $t=0,\ldots,N$.
The index $t$ denotes time, and the system evolves according to linear dynamics constraint (\ref{eq:time_split_objective_dyn}) where $c_t \in \reals^n$ is considered to be a
known disturbance.
Here, $N$ is the prediction horizon and $Q \in \reals^{n\times n}, R \in \reals^{m\times m}$ are symmetric matrices.
The stage cost functions in (\ref{eq:time_split_objective_cost}) are convex quadratic with $Q\succeq 0$ and $R\succ 0$.
The stage-wise state-input pairs are constrained to reside within polyhedra (\ref{eq:time_split_objective_constr}) denoted by $\mathcal{X}_t$ and $\mathcal{U}_t$, respectively.
These are constraint sets defined by linear inequalities that involve states and inputs at the same sample time index.  \\

\begin{figure}[!htb]
      \begin{center}
         \includegraphics[scale=0.15]{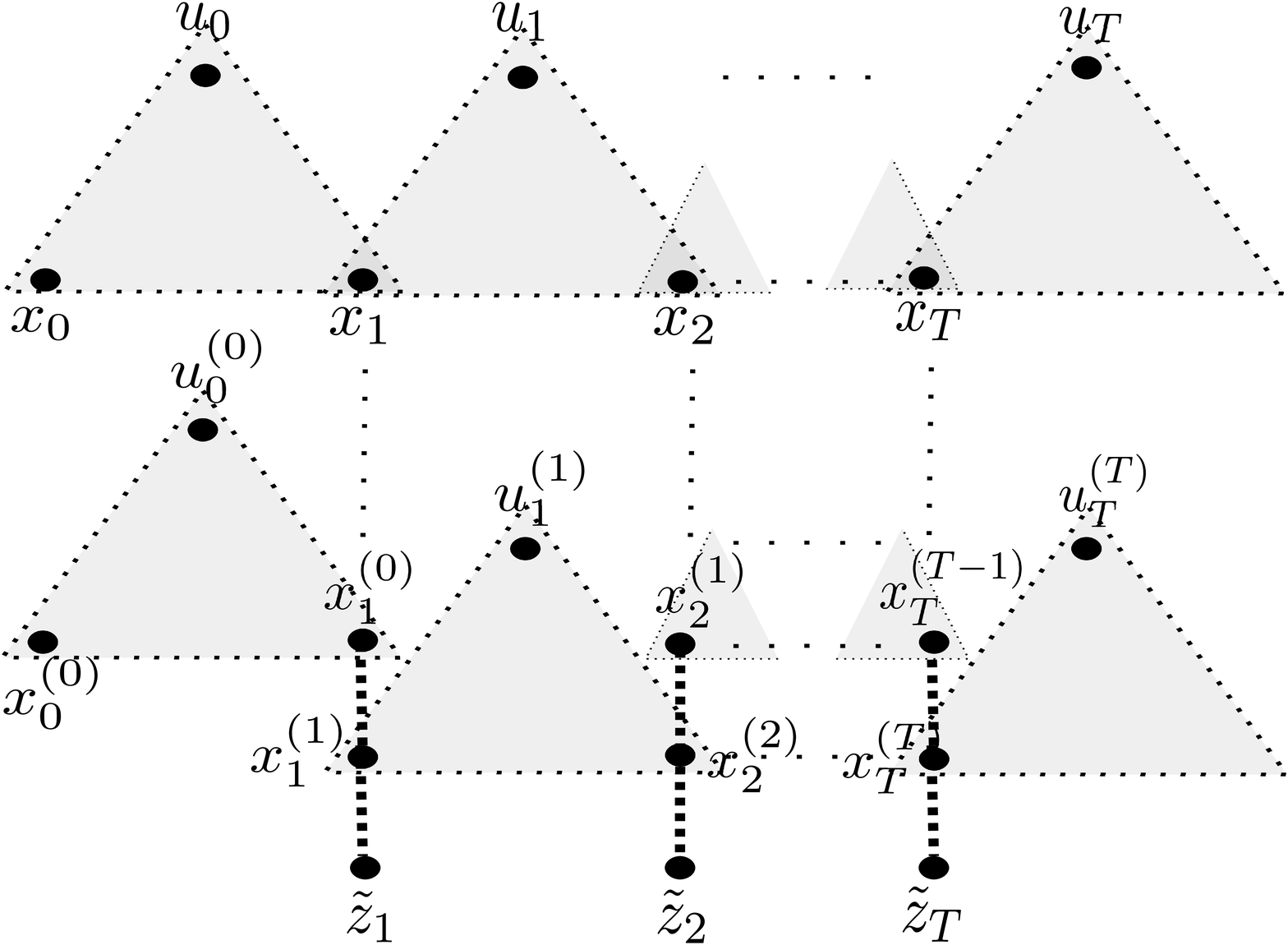}
	 \caption{The idea of the time splitting algorithm. Whenever the dynamics introduce a coupling
           of variables, they are decoupled via a slack variable $\tilde{z}_t=x_t^t=x_t^{(t-1)}$.}
         \label{time_split}
      \end{center}
\end{figure}

Motivated by the principles of operator splitting methods (see \cite{BPCPE2011} for details and relevant references),
we propose to split the problem (\ref{eq:time_split_objective}) into $N+1$ smaller stage-wise subproblems that can be solved in parallel.
This requires breaking the coupling that appears due to the dynamics.
We introduce a copy of each variable that couples the dynamics equations in order to allow such a splitting into subproblems, and subsequently
impose a consensus constraint on the associated complicating variables and their copies. This leads to the following equivalent formulation of
(\ref{eq:time_split_objective}), where the complicating variables that are used to perform the splitting are clearly highlighted:
	\begin{subequations}	\label{eq:time_split_objective2}
	\begin{eqnarray}
	  \mbox{minimize} && \frac{1}{2}\sum_{t=0}^N \Big(x_t^{(t)T} Q x_t^{(t)} + u_t^{(t)T} R u_t^{(t)}\Big) \\
	  \mbox{subject to} && (x_t^{(t)},u_t^{(t)}) \in \mathcal{X}_t \times \mathcal{U}_t, \quad t = 0,\ldots,N \\
		&& x_{t+1}^{(t)} = A_t x_t^{(t)}+B_t u_t^{(t)}+c_t,\\
		&& t = 0,\ldots,N-1 \nonumber \\
		&& x_0^{(0)} = x_\mathrm{init} \\
		&& \tilde{z}_{t+1} = x_{t+1}^{(t)} \\
		&& \tilde{z}_{t+1} = x_{t+1}^{(t+1)}, \quad t = 0,\ldots,N-1,
	\end{eqnarray}
	\end{subequations}
where the subscript $t$ of the decision variables $x_t^{(t)},u_t^{(t)}$ indicates the time index and the superscript $(t)$ denotes the group or subproblem
where the variable belongs to. Hence, each subproblem contains three variables, the current state and input as well as a prediction of the
state for the next time instant. 
The introduced complicating variable $\tilde{z}_t$ acts as a `global' variable that brings the local copies $x_t^{(t-1)}$
and $x_t^{(t)}$ in agreement, i.e., $\tilde{z}_t=x_t^{(t-1)}=x_t^{(t)}$. The time-splitting idea is graphically depicted in
Figure~\ref{time_split}.

\subsection{The time-splitting algorithm} \label{subsec:timesplitting}

In order to use a more compact formulation, we will denote the decision variables in (\ref{eq:time_split_objective2}) corresponding to each subproblem using
\begin{equation} \label{eq:time_split_augmented_vec}
 \tilde{x}_t = (x_t^{(t)}, u_t^{(t)}, x_{t+1}^{(t)}),
\end{equation}
where $\tilde x_t \in \reals^{2n+m}$. We also introduce dual variables to deal with the consensus equality constraints:
\begin{itemize}
 \item $\tilde{w}_t$ associated with $x_t^{(t)}=\tilde{z}_t, \quad t=1,\ldots,N$ \text{and} \\
 \item $\tilde{v}_t$ associated with $x_t^{(t-1)}=\tilde{z}_t, \quad t=1,\ldots,N$.
\end{itemize}
In order to rewrite the finite-time optimal control problem in a more compact form, we define the following matrices:
\begin{subequations} \label{time_split_matrices}
\begin{eqnarray}
 P_t &=& \left[\begin{array}{ccc} Q_t & 0 & 0 \\ 0 & R_t & 0 \\ 0 & 0& 0\end{array}\right] \in \reals^{(2n+m)\times (2n+m)},  \\
 F_t &=& \left[\begin{array}{ccc} -A_t & -B_t & I\end{array}\right] \in \reals^{n\times (2n+m)}, \\
 G_0 &=& \left[\begin{array}{ccc} I & 0 & 0\end{array}\right] \in \reals^{n\times (2n+m)}, \\
 G_1 &=& \left[\begin{array}{ccc} 0 & 0 & I\end{array}\right] \in \reals^{n\times (2n+m)}.
\end{eqnarray}
\end{subequations}


We use the ALternating Direction Method of Multipliers (ADMM) \cite{GlM:75,GaM:76} in order to arrive at a solution approach that is amenable to parallel implementation.
The updates involved in the ADMM algorithm include forming the augmented Lagrangian of the problem and minimizing over the primal variables $\tilde{x}_t, \; t=0,\ldots,N$ and
$\tilde{z}_t, \; t=1,\ldots,N$, followed by updating the dual variables $\tilde{w}_t$ and $\tilde{v}_t, \; t=1,\ldots,N$.
The three main steps of the algorithm are performed in an iterative fashion and are described next in detail. We use $k$ to denote the algorithm's loop counter. The termination criterion based on primal and dual tolerances are provided and the
analytic derivation of the formulas in each step is presented in the Appendix. \\


\textbf{Step 1: Solving $N+1$ QP subproblems for the primal variables in $\tilde{x}_t$} \\
Minimization of the augmented Lagrangian over the primal variables $\tilde{x}_t$ results in $N+1$ stage-wise
quadratic programs (QPs):
\begin{itemize}

\item For the subproblem associated with the time instant $t=0$, we need to solve
 \begin{equation}  \label{eq:time_split_x0}
 \begin{array}{llll}
    & \mbox{minimize}
    && (1/2)\tilde x_0^TP_0 \tilde x_0 -\rho \tilde v_1^{kT}(G_1\tilde x_0-\tilde z_1^k)  \\
    &&&     +(\rho/2)\|G_1\tilde x_0-\tilde z_1^k\|_2^2 \\
    & \mbox{subject to}
    & & \tilde x_0 \in \mathcal{C}_0 \\
    &&& G_0\tilde x_0 = x_\mathrm{init} \\
    &&& F_0\tilde x_0 = c_0
 \end{array}
 \end{equation}
 with variable $\tilde x_0$.

\item Similarly, we need to solve the following QPs for all the other groups of variables $\tilde x_t, \; t=1,\ldots,N-1$:
 \begin{equation} \label{eq:time_split_xt}
 \begin{array}{llll}
    & \mbox{minimize}
    && (1/2)\tilde x_t^TP_t \tilde x_t -\rho \tilde v_{t+1}^{kT}(G_1\tilde x_t-\tilde z_{t+1}^k)  \\
    &&&   +(\rho/2)\|G_1\tilde x_t-\tilde z_{t+1}^k\|_2^2 \\
    &&& -\rho \tilde w_t^{kT}(G_0\tilde x_t-\tilde z_t^k) \\
    &&&   +(\rho/2)\|G_0\tilde x_t-\tilde z_t^k\|_2^2\\
    & \mbox{subject to}
    &&  \tilde x_t \in \mathcal{C}_t \\
    &&& F_t\tilde x_t = c_t
 \end{array}
 \end{equation}

 \item For the subproblem associated with the final time instant $t=N$, we need to solve
 \begin{equation} \label{eq:time_split_xT}
 \begin{array}{llll}
    & \mbox{minimize}
    && (1/2)\tilde x_N^TP_N \tilde x_N -\rho \tilde w_N^{kT}(G_0\tilde x_N-\tilde z_N^k)  \\
    &&&   +(\rho/2)\|G_0\tilde x_N-\tilde z_N^k\|_2^2\ \\
    & \mbox{subject to}
    & & \tilde x_N \in \mathcal{C}_N
 \end{array}
 \end{equation}
 with variable $\tilde x_N$.
\end{itemize}
The polyhedral sets $\mathcal{C}_t, \; t=0,\ldots,N$ are defined as
\begin{equation}
\mathcal{C}_t =\mathcal{X}_t \times \mathcal{U}_t \times \mathcal{X}_{t+1} \subseteq \reals^{2n+m}
\end{equation}
and the variable $\rho > 0$ is a parameter of the algorithm.


\begin{remark}
Notice that for the time instant $N$ the decision variables of the QP actually simplify to $\tilde x_N = x_N$ and
$\mathcal{C}_N = \mathcal{X}_N$, but we keep the same notation for simplicity (and without loss of generality).
\end{remark}

\textbf{Step 2: Averaging} \\
The update of the `global' primal variables $\tilde{z}_t, \; t=1,\ldots,N$ is derived from a simple quadratic minimization
problem, the solution of which turns out to be an average of the predicted ($x_t^{(t-1)}$) and current ($x_t^{(t)}$) state
\begin{equation} \label{step2_time_split}
 \tilde z_t^{k+1} = \frac{G_0\tilde x_t^{k+1}+G_1\tilde x_{t-1}^{k+1}}{2}, \quad t=1,\ldots,N.
\end{equation}
This intuitively makes sense, since the global variable can be obtained by collecting the local (primal) ones and computing the best estimate based on their values.

\textbf{Step 3: Dual update} \\
The dual updates can be expressed as
\begin{subequations} \label{step3_time_split}
\begin{eqnarray}
 \tilde w_t^{k+1} &=& \tilde w_t^k - G_0\tilde x_t^{k+1} + \tilde z_t^{k+1},\\
 \tilde v_t^{k+1} &=& \tilde v_t^k - G_1\tilde x_{t-1}^{k+1} + \tilde z_t^{k+1}, \quad t=1,\ldots,N.
\end{eqnarray}
\end{subequations} \\

\textbf{Termination criterion} \\
The algorithm terminates when a set of \emph{primal} and \emph{dual} residuals are bounded by a specified threshold
(primal and dual tolerances); see \cite[\S 3.2]{BPCPE2011}. The primal and dual residuals for the time-splitting algorithm are
defined as
\begin{equation}
r^k = A_{\mathrm{res}}x_\mathrm{pri}^k + B_{\mathrm{res}}z_\mathrm{pri}^k,
\end{equation}
and
\begin{equation}
s^k = -\rho A_{\mathrm{res}}^TB_{\mathrm{res}}(z_\mathrm{pri}^k-z_\mathrm{pri}^{k-1}),
\end{equation}
respectively.

The termination criterion is activated when
\[
\|r^k\|_2 \leq \epsilon^\mathrm{pri}, \quad \|s^k\|_2 \leq \epsilon^\mathrm{dual},
\]
where the tolerances $\epsilon^\mathrm{pri}$ and $\epsilon^\mathrm{dual}$
are defined as follows:
\begin{subequations} \label{tols-time_split}
\begin{eqnarray}
\epsilon^\mathrm{pri} &=&  \epsilon^\mathrm{abs} \sqrt{N2n}  \\
 &~& +\epsilon^\mathrm{rel} \max\{ \|A_{\mathrm{res}}x_\mathrm{pri}\|_2,\|B_{\mathrm{res}}z_\mathrm{pri}\|_2,\|c_{\mathrm{res}}\|_2 \} \\
\epsilon^\mathrm{dual} &=&  \epsilon^\mathrm{abs} \sqrt{(2n+m)(N+1)} + \epsilon^\mathrm{rel} \|A_{\mathrm{res}}^Tv_\mathrm{dual}\|_2
\end{eqnarray}
\end{subequations}
and we defined the vectors
\[
x_\mathrm{pri} = \left[\begin{array}{c} \tilde x_0 \\ \tilde x_1 \\ \vdots \\ \tilde x_N \end{array} \right], \qquad
z_\mathrm{pri} = \left[\begin{array}{c} \tilde z_1 \\ \tilde z_2 \\ \vdots \\ \tilde z_N \end{array} \right], \qquad
v_\mathrm{dual} = \left[\begin{array}{c} \tilde v_1 \\ \tilde w_1 \\ \tilde v_2 \\ \vdots \\ \tilde v_N \\
 \tilde w_N \end{array} \right].
\]
The residual matrices $A_{\mathrm{res}}$ and $B_{\mathrm{res}}$ and the residual vector $c_{\mathrm{res}}$ are
\[
A_{\mathrm{res}} = \left[\begin{array}{cccccc}
G_1 & 0 & 0 & \cdots & 0 & 0 \\
0 & G_0 & 0 & \cdots & 0 & 0 \\
0 & G_1 & 0 & \cdots & 0 & 0 \\
0 & 0 & G_0 & \cdots & 0 & 0 \\
0 & 0 & G_1 & \cdots & 0 & 0 \\
\vdots & \vdots & \vdots & \ddots & \vdots & \vdots \\
0 & 0 & 0 & \cdots & G_0 & 0 \\
0 & 0 & 0 & \cdots & G_1 & 0 \\
0 & 0 & 0 & \cdots & 0 & G_0
\end{array} \right] \in \reals^{T2n \times (N+1)(2n+m)},
\]
\[
B_{\mathrm{res}} = \left[\begin{array}{cccc}
-I & 0 & \cdots & 0 \\
-I & 0 & \cdots & 0 \\
0 & -I & \cdots & 0 \\
0 & -I & \cdots & 0 \\
\vdots & \vdots & \ddots & \vdots \\
0 & 0 & \cdots & -I \\
0 & 0 & \cdots & -I \\
\end{array} \right] \in \reals^{N2n \times Nn},
\]
and $c_{\mathrm{res}}\in \reals^{N2n}$ is zero.\\

The three update steps described above are fully parallelizable at each iteration $k$. Assuming $N+1$ processors are available, then processor $\Pi_t$ would need
to execute the following actions for $t=0,\ldots,N-1$:
\begin{enumerate}
 \item Receive the estimate $\tilde z_{t+1}^k$ $\tilde v_{t+1}^k$ and from neighboring processor $\Pi_{t+1}$ (\ref{eq:time_split_xt}).
 \item Compute $\tilde x_t^{k+1}$.
 \item Receive the estimate $\tilde x_{t-1}^{k+1}$ from neighboring processor $\Pi_{t-1}$ and compute $\tilde z_t^{k+1}$ (\ref{step2_time_split}).
 \item Compute $\tilde w_t^{k+1}$ and $\tilde v_t^{k+1}$ (\ref{step3_time_split}).
 \item Communicate $\tilde x_t^{k+1}$ to processors $\Pi_{t-1}$ and $\Pi_{t+1}$, and $\tilde z_t^{k+1}$, $\tilde v_t^{k+1}$ to processor $\Pi_{t-1}$.
\end{enumerate}
The above scheme suggests that each processor $\Pi_t, \; t=1,\ldots,N-1$ interacts with the two neighboring processors
$\Pi_{t-1}$ and $\Pi_{t+1}$. Processors $\Pi_0$ and $\Pi_N$ communicate only with processors $\Pi_1$ and $\Pi_{N-1}$,
respectively. After updating all variables, a gather operation follows in order to compute the residuals and check the termination
criterion.

\section{THREE-SET SPLITTING QP SOLVER} \label{sec:threeset}

\subsection{Motivation}

The time-splitting algorithm presented in the previous section decomposes the centralized finite-time optimal control problem
so that it can be solved using multiple parallel processors. However, the updates for the primal variables $\tilde{x}_t, \; t=0,\ldots,N$ given in
(\ref{eq:time_split_x0}), (\ref{eq:time_split_xt}) and (\ref{eq:time_split_xT}) involve solving a QP at each iteration of the
algorithm. Even though several fast interior point solvers exist for this purpose (see e.g., \cite{MB2012}), these are mostly suitable for
only a limited number of variables. Although recently more computationally efficient schemes that scale better with the problem size have been developed,
they are restricted to cases where simple box constraints are considered \cite{WB2010,OSB2012,KF2012}.
In order to achieve fast computations in an embedded control environment, other generic solution methods would be preferred.

In this section we propose an alternative scheme for solving the QPs with general polyhedral
constraints that arise in the previous section. We propose to perform yet another type of splitting approach, which splits the state-input
variables of the QP in three sets. One set involves the variables that appear in the
objective function, another set includes those that appear in the dynamics equality constraints,
and the last set contains variables from the inequality constraints. In this way, we solve three simpler subproblems instead of the
single general QP. Since several variables are shared among the subproblems, their solutions must be in consensus again to ensure consistency. \\

An important element of the proposed method is the introduction of an extra slack variable,
which allows to get an analytic solution for the subproblem associated with the inequality constraints.
Using this variable, the projection on any polyhedral set can be practically rewritten as a projection onto the nonnegative orthant.
Besides this feature, the proposed splitting exploits structure in the resulting matrices and thus modern numerical linear algebra
methods can be employed for speeding up the computations. \\

\subsection{Problem setup}

We consider a QP of the form	
 \begin{equation} \label{eq:quadratic_objective}
 \begin{array}{ll}
    \mbox{minimize} & \frac{1}{2}x^TMx + q^Tx + r\\
    \mbox{subject to} & Ax = b \\
    \quad & Hx \preceq h,
 \end{array}
 \end{equation}
with decision variable $x \in \reals^n$, where $M \in S_+^n$, $q \in \reals^n$, $r \in \reals$, $A \in \reals^{m\times n}$,
$b \in \reals^m$, $H \in \reals^{p\times n}$, $h \in \reals^p$ and $S_+^n$ is the cone
of positive semidefinite matrices of dimension $n$.

In order to apply a variable splitting idea for this problem we first replicate all variables appearing in (\ref{eq:quadratic_objective})
three times, introducing three different sets for which we must ensure consensus. Furthermore, we use a slack variable to remove
the polyhedral constraint and transform it into a projection operation onto the nonnegative orthant.
We define the following sets of variables:
\begin{itemize}
 \item First set - objective: $x^{\I}$ \\
 \item Second set - equality constraints: $x^{\II}$ \\
 \item Third set - inequality constraints: $x^{\III}$ \\
 \item First global variable: $z=x^{\I}=x^{\II}=x^{\III}$ \\
 \item Second global variable: $y=h-Hx^{(\III)}$ \\
 \item First set of dual variables: $\tilde z^{\I},\tilde z^{\II},\tilde z^{\III}$ associated with
            $z=x^{\I}=x^{\II}=x^{\III}$, respectively.\\
 \item Second set of dual variables: $\tilde y$ associated with $y=h-Hx^{\III}$
\end{itemize}

Using the above sets of variables problem (\ref{eq:quadratic_objective}) can be restated in the equivalent form
 \begin{subequations} \label{eq:3sets_split_objective}
 \begin{eqnarray}
    \mbox{minimize} && \frac{1}{2}x^{\I T}Mx^{\I} + q^Tx^{\I} + r\\
    \mbox{subject to} && Ax^{\II} = b \\
    \quad && y = h - Hx^{\III}, \quad y \succeq 0 \\
    \quad && x^{\I} = x^{\II} = x^{\III} = z,
 \end{eqnarray}
 \end{subequations}
with variables $x^{\I}, x^{\II}, x^{\III}, z \in \reals^n, y \in \reals^p$. The dual variables
are of dimensions $\tilde{z} \in \reals^n, \tilde{y} \in \reals^p$.

\subsection{The proposed three-set splitting algorithm}

The proposed algorithm consists of iterative updates to the three ADMM steps similarly to the case of the time-splitting approach in Section~\ref{subsec:timesplitting},
namely one for the (local) primal variables $x^{\I}, x^{\II}, x^{\III}$, one for the (global) primal variables $z$ and $y$ and one for the dual variables
$\tilde z^{\I},\tilde z^{\II},\tilde z^{\III}$ and $\tilde y$. We provide the algorithm's steps below along with some clarifying
comments. The analytic derivations are presented in the Appendix. \\


\textbf{Step 1: Solving three subproblems for the primal variables $x^{\I},x^{\II},x^{\III}$} \\
In all three cases we have to solve simple, unconstrained QPs. The updates are:

 \begin{equation}\label{eq:3_sets_x1}
  \left(x^{\I}\right)^{k+1} = (M+\rho I)^{-1}(\rho(z^k +\left(\tilde z^{\I}\right)^k )-q)
 \end{equation}
 \begin{equation}\label{eq:3_sets_x2}
  \left[\begin{array}{cc} \rho I &  A^T \\  A & 0\end{array}\right]
  \left[\begin{array}{c}\left(x^{\II}\right)^{k+1} \\ \nu\end{array}\right] =
  \left[\begin{array}{c}\rho(z^k+\left(\tilde z^{\II}\right)^k) \\ b\end{array}\right],
 \end{equation}
where $\nu$ is the dual variable associated with the equality constraint $Ax^{\II}=b$, and
 \begin{eqnarray}
   \left(x^{\III}\right)^{k+1} \nonumber &=& (H^TH+ I)^{-1}\Big(H^T( h- \tilde y^k-y^k)+ \\
                              && z^k+\left(\tilde z^{\III}\right)^k\Big){\label{eq:3_sets_x3}}.
 \end{eqnarray}
The matrices $M+\rho I$ and $H^TH+I$ are symmetric positive definite due to the regularization terms. This means that,
instead of directly inverting the matrices, we can save computational effort by taking the Cholesky
factorization, i.e., write the matrix as a product of a lower triangular matrix and its transpose (see, e.g., \cite{GVL1966}).
Furthermore, the matrix $\left[\begin{array}{cc} \rho I & A^T \\ A & 0\end{array}\right]$ is a KKT matrix and
$\rho I \succ 0$. Hence we can exploit its structure and use block elimination to solve the KKT system
(see \cite[App. C]{BV2004}). The resulting matrices can be pre-factorized and then used in every solve step. The right-hand sides are the only parts that change in
the update loop.

\textbf{Step 2: Averaging and projection} \\
The update for $z$ is
\begin{equation}
 z^{k+1} = \frac{1}{3}\sum_{i=\I}^{\III} \left(x^{i}\right)^{k+1},
\end{equation}
while for $y$ it is
\begin{equation}
  y^{k+1} = \left(h-H\left(x_t^{\III}\right)^{k+1}-\tilde y^k\right)_+.
\end{equation}
The $z$-update is an averaging over the three sets of the primal variables $x^{\I}, x^{\II}, x^{\III}$, while the
$y$-update is the solution of a proximal minimization problem (see Appendix), resulting in
a projection onto the nonnegative orthant, denoted by $(\cdot)_+$.

\textbf{Step 3: Dual update} \\
The update for the dual variables $\tilde z^{i}$ is
\begin{equation}
 \left(\tilde z^{i}\right)^{k+1} = \left(\tilde z^{i}\right)^k - \left(x^{i}\right)^{k+1} + z^{k+1}, \quad i=\I, \II, \III.
\end{equation}

Similarly, the update for $\tilde y$ is
\begin{equation}
 \tilde y^{k+1} = \tilde y^k + y^{k+1} - h + H\left(x^{\III}\right)^{k+1}.
\end{equation}

\textbf{Termination criterion} \\
The primal and dual tolerances $\epsilon^\mathrm{pri}$ and $\epsilon^\mathrm{dual}$ are given by
\begin{equation}\label{tols-3sets}
\begin{array}{rcl}
\epsilon^\mathrm{pri} &=&  \epsilon^\mathrm{abs} \sqrt{3n+p}+ \epsilon^\mathrm{rel}
\max\{  \|A_{\mathrm{res}}(x^{\I},x^{\II},x^{\III})\|_2, \\
 && \|B_{\mathrm{res}}(z,y)\|_2, \|c_{\mathrm{res}}\|_2 \} \\
\epsilon^\mathrm{dual} &=&  \epsilon^\mathrm{abs} \sqrt{3n} + \epsilon^\mathrm{rel}
\|A_{\mathrm{res}}^T(\tilde z^{\I},\tilde z^{\II},\tilde z^{\III},\tilde y)\|_2.  \end{array},
\end{equation}
The residual matrices $A_{\mathrm{res}}$ and $B_{\mathrm{res}}$ and the residual vector $c_{\mathrm{res}}$ are
\begin{equation*}
A_{\mathrm{res}} = \left[\begin{array}{ccc}
I & 0 & 0 \\
0 & I & 0 \\
0 & 0 & I \\
0 & 0 & H
\end{array} \right], \quad
B_{\mathrm{res}} = \left[\begin{array}{cc}
-I & 0 \\
-I & 0 \\
-I & 0 \\
 0 & I
\end{array} \right],
\end{equation*}
and
\[
 c_{\mathrm{res}} = (0,0,0,h),
\]
where $A_{\mathrm{res}} \in \reals^{(3n+p) \times 3n}$, $B_{\mathrm{res}} \in \reals^{(3n+p) \times (n+p)}$
and $c_{\mathrm{res}} \in \reals^{3n+p}$.

\begin{remark}
For the QP corresponding to the last sample time $t=N$ (\ref{eq:time_split_xT}),
the algorithm simplifies to splitting into two sets (objective and inequality constraints), since there are
no dynamics equality constraints. The updates and residuals follow directly from the more generic case presented above.
\end{remark}

\section{HIERARCHICAL TIME-SPLITTING OPTIMAL CONTROL}   \label{sec:hierarchical}

\begin{figure*}[!htb]
      \begin{center}
         \includegraphics[scale=0.75]{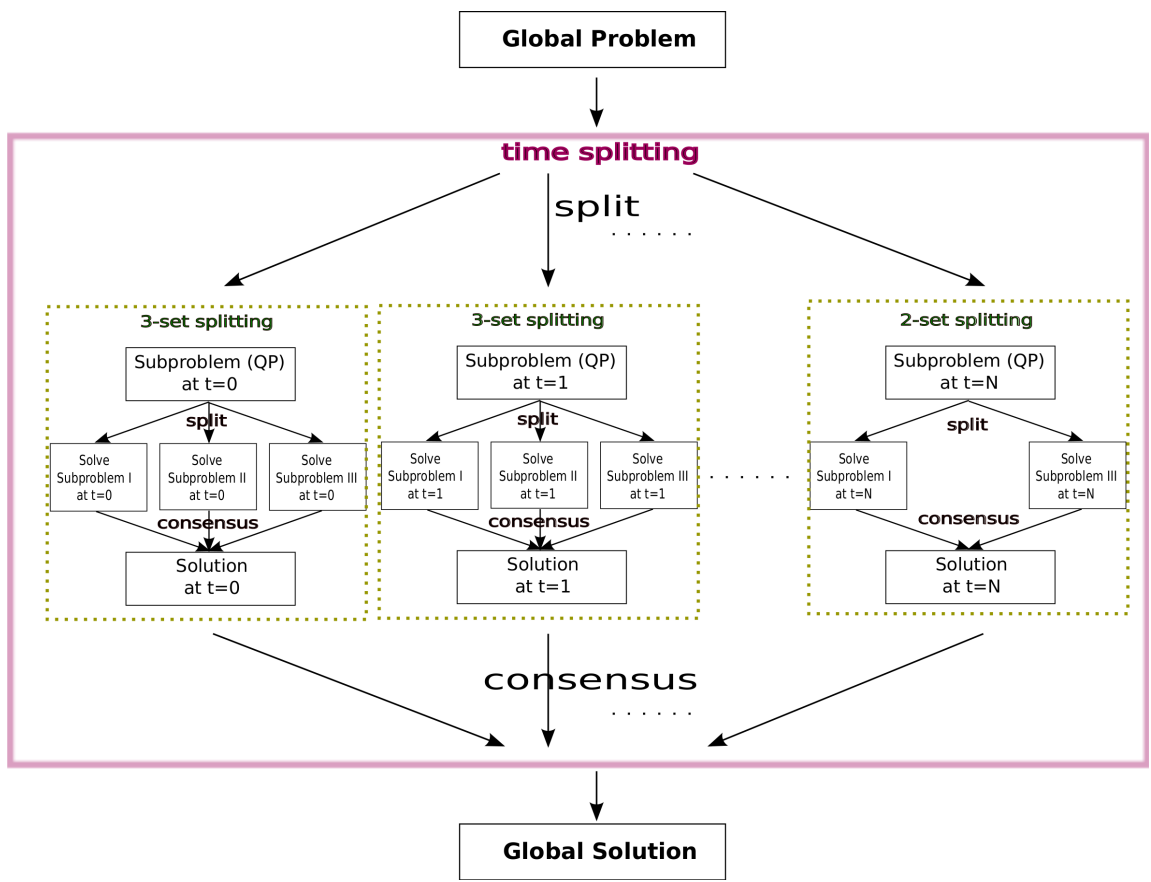}
	 \caption{Structural representation of the hierarchical time-splitting solution approach to the finite-time optimal control problem.
     The outer subproblems at $t=0, 1, \ldots, N$ refer to (\ref{eq:time_split_x0}), (\ref{eq:time_split_xt}) and (\ref{eq:time_split_xT}), respectively.
     The solutions of the nested subproblems enumerated by $\I$, $\II$ and $\III$ correspond to (\ref{eq:3_sets_x1}), (\ref{eq:3_sets_x2}) and (\ref{eq:3_sets_x3}).
     Notice that only 2 subproblems have to be solved for the last time instant $t=N$. }
         \label{hierarch_time_split_ctrl}
      \end{center}
\end{figure*}

It is a natural idea to combine the two splitting algorithms (time-splitting and three-set splitting) that were introduced in the preceding two sections
in order to speed up the solution of the finite-time optimal control problem (\ref{eq:time_split_objective}).
This can be accomplished via a nested decomposition scheme, where we employ the
\emph{three-set splitting} algorithm to solve the QPs (\ref{eq:time_split_x0}), (\ref{eq:time_split_xt})
and (\ref{eq:time_split_xT}) appearing in Step~1 of the \emph{time-splitting} algorithm. The idea is graphically depicted
in Figure~\ref{hierarch_time_split_ctrl}. 
  \\

If we rewrite the generalized inequality constraints appearing in the problem formulation as
\begin{equation}
 (x_t^{(t)},u_t^{(t)}) \in \mathcal{C}_t \; \Leftrightarrow \; H_t \tilde x_t \preceq h_t,
                 \quad t = 0,\ldots,N,
\end{equation}
then the QPs can be written in the form described by (\ref{eq:quadratic_objective}), where we consider the following
relations:
\begin{itemize}
 \item For $\tilde x_0$, Eq. $(\ref{eq:time_split_x0})$: \\
       \begin{eqnarray*}
	  M &=& P_0+\rho G_1^TG_1, \\
	  q &=& -\rho G_1^T(\tilde z_1 + \tilde v_1), \\
	  r &=& (\rho/2)\tilde z_1^T\tilde z_1 + \rho \tilde v_1^T\tilde z_1, \\
	  A &=& (G_0, F_0), \\
          b &=& (x_\mathrm{init}, c_0), \\
          H &=& \left[\begin{array}{cc} H_0 & 0 \end{array}\right], \\
          h &=& h_0.
       \end{eqnarray*}

 \item For $\tilde x_t, \; t=1,\ldots,N-1$, Eq. $(\ref{eq:time_split_xt})$: \\
       \begin{eqnarray*}
	  M &=& P_t+\rho (G_0^TG_0+G_1^TG_1), \\
	  q &=& -\rho (G_0^T(\tilde z_t + \tilde w_t) + G_1^T(\tilde z_{t+1} + \tilde v_{t+1})), \\
	  r &=& (\rho/2)(\tilde z_{t+1}^T\tilde z_{t+1} + \tilde z_t^T\tilde z_t) +
                         \rho (\tilde v_{t+1}^T\tilde z_{t+1} + \tilde w_t^T\tilde z_t), \\
	  A &=& F_t, \\
          b &=& c_t, \\
          H &=& \left[\begin{array}{cc} H_t & 0 \end{array}\right], \\
          h &=& h_t.
       \end{eqnarray*}

 \item For $\tilde x_N$, Eq. $(\ref{eq:time_split_xT})$: \\
       \begin{eqnarray*}
	  M &=& P_N+\rho G_0^TG_0, \\
	  q &=& -\rho G_0^T(\tilde z_N + \tilde w_N), \\
	  r &=& (\rho/2)\tilde z_N^T\tilde z_N + \rho \tilde w_N^T\tilde z_N, \\
	  A &=& 0, \\
          b &=& 0, \\
          H &=& \left[\begin{array}{cc} H_N & 0 \end{array}\right], \\
          h &=& h_N.
       \end{eqnarray*}
\end{itemize}

For each iteration of the time-splitting algorithm, the
three-set splitting algorithm runs in an inner loop until it
converges. The quality of this convergence, i.e., the choice
of the primal and dual tolerances of the inner loop (\ref{tols-3sets})
will affect the quality of the global solution. 

A method that enables substantial speedup of the algorithm is warm starting.
Since the three-set splitting algorithm will run for every iteration of the
time-splitting algorithm, we can warm-start each QP with its previous solution.
In this way, we can achieve a significant reduction in the number of iterations
needed for the convergence of the inner loop.

\section{NUMERICAL RESULTS}  \label{sec:results}

We consider three, randomly generated, numerical examples to illustrate the
performance of the algorithm. The examples vary in terms of the number of decision variables involved.
The systems considered are linear and time-invariant. We impose constraints on the difference between
two consecutive states at each time instant of the form
\[
 x_{t,i}-x_{t,i-1} \leq dx,
\]
where $x_{ti} \in \reals, \; i=1, \ldots,n, \;, t=0,\ldots,N$ and box constraints on the inputs, i.e.,
\[
 \|u_t\|_{\infty} \leq u_\mathrm{max}, \; t=0,\ldots,N-1.
\]
By adjusting the level of the disturbance $c_t$ in several time instances, we ensure activation of the constraints along the horizon.

For the simulations, we used an Intel Core i7 processor running at 1.7 GHz.
We compared a C-implementation of our algorithm with using CVX \cite{GBY2006}, a parser-solver
that uses SDPT3 \cite{TTT1999}. For our method, the tolerances for both the outer
and inner algorithms are set as $\epsilon^\mathrm{pri}=10^{-4}$ and
$\epsilon^\mathrm{dual}=10^{-3}$. The parameter $\rho$ was set after some simple tuning.
The linear systems appearing in (\ref{eq:3_sets_x1}), (\ref{eq:3_sets_x2}) and (\ref{eq:3_sets_x3})
were solved by first factorizing the matrices off-line, using
Tim Davis's sparse package \cite{ldl,amd,davis_book} (see also \cite{OSB2012}). The
finite-time optimal control problem was solved only once and all the
primal and dual variables were initialized at zero. However,
the inner algorithm was warm-started at every iteration of
the outer algorithm to the values acquired from the previous
iteration. No relaxation or any other variance of the iterations
was used. The numerical results are summarized in Table~\ref{Table},
where the computation times are reported in ms.

\begin{table}
\normalsize
\begin{center}
\scalebox{0.85}{
\begin{tabular}{|c||c|c|c|}
\hline
 & \textbf{small} & \textbf{medium}
& \textbf{large}\\
\hline\hline
states $n$ & 10 & 20  & 50\\
\hline
inputs $m$ & 10 & 10 & 40 \\
\hline
horizon length $N$ & 10 & 30 & 60 \\
\hline
total variables & 220  &900 & 5400\\
\hline
$\rho$ & 15 & 25 & 50 \\
\hline
active box constraints & 5  & 6 & 20\\
\hline
active inequality constraints & 2  & 2 & 4\\
\hline\hline
CVX solve time & 2430 & 3529 & 19420\\
\hline\hline
factorization time & 3.18 & 9.5  & 30\\
\hline\hline
&\multicolumn{3}{|c|}{Tolerance $10^{-4}$} \\
\hline
3-set (average) iterations & 21.80  & 17  & 15.95 \\
\hline
3-set (average) solve time & 0.75 & 1.38 & 9.15\\
\hline
time-split. iterations & 250 & 241 & 389\\
\hline
time-split. solve time (single thread) & 1880 &  10023  & 215888\\
\hline
time-split. solve time ($N$ threads$^\ast$) & 188 &  334.1  & 3598\\
\hline\hline
&\multicolumn{3}{|c|}{Tolerance $10^{-3}$} \\
\hline
3-set (average) iterations & 13.14  & 13.27  &12.32 \\
\hline
3-set (average) solve time & 0.49 & 1.18 & 7.34\\
\hline
time-split. iterations & 156 & 128 & 224\\
\hline
time-split. solve time (single thread) & 780 &  4304  & 99525\\
\hline
time-split. solve time ($N$ threads$^\ast$) & 78 &  143.47  & 1659\\
\hline \hline
\multicolumn{4}{c}{$^\ast$ estimated parallel computation times}
\end{tabular}
}
\end{center}
\caption{Hierarchical time-splitting optimal control: computational time results for different size problems.}
\label{Table}
\end{table}

We can observe that, in the case of the small system, even
when solving the problem on a single thread, the computation times
are smaller than those of CVX. As the problem scales,
the computations have to be parallelized in order to gain a significant
advantage. More specifically, we expect the following
speedup factors: 13 and 31 times faster in the case of the small problem
(for the corresponding tolerances set to $10^{-4}$ and $10^{-3}$
respectively). For the medium-sized problem the speedups
are by a factor of 10.5 and 24.6, and a factor of 5.4 and 11.7 for the
large-scale problem, respectively.

In addition, we could observe that the factorization times are
negligible in all cases, since the matrices being factorized
are not large. Concerning the three-set splitting algorithm, only the average
computational times are indicated over all iterations required to solve the problem.

\section{CONCLUSIONS} \label{sec:conclusions}


In this paper, we proposed an algorithm that solves
a centralized convex finite-time optimal control problem
making use of operator splitting methods, and, more
specifically, the Alternating Direction Method of Multipliers.
The initial problem is split into as many subproblems as the
horizon length, which then can be solved in parallel.

The resulting algorithm is composed of three steps,
including one where several QPs have to be solved. In this respect, we
proposed another method, based again on operator splitting,
that is applicable to QPs of any size, involving polyhedral
constraints. This algorithm exploits the structure of the
problem, leading to fast solutions.

The combination of the proposed algorithms results in a
nested decomposition scheme for solving the aforementioned finite-time
optimal control problems over several parallel processors.

Our numerical experiments suggest that the proposed hierarchical decomposition
approach provides significant speed-up in computational time required for medium
accuracy solutions for the class of problems considered. In our future work
we intend to perform an even more extensive comparison with very recent tailor-made
computational tools, and implement the algorithm on a parallel computing platform
to obtain more accurate and representative computational time measurements.



\section*{APPENDIX}

\subsection*{Derivation of the time-splitting algorithm updates}

We solve the relaxed version of the convex optimization problem (\ref{eq:time_split_objective2}) by formulating the
augmented Lagrangian with respect to the additional equality constraints, using
the matrices and vectors defined in (\ref{eq:time_split_augmented_vec}) and (\ref{time_split_matrices}).
The augmented Lagrangian can be written as

\begin{eqnarray}
&& L_{\rho}(\tilde x_t,x_\mathrm{init},\tilde w_t,\tilde y_t,\tilde z_t,\nu_t) := \label{eq:AugLag_time_split} \\
&& ~~~~~~ \nonumber \sum_{t=0}^N \left( \frac{1}{2}\tilde x_t^TP\tilde x_t + I_{\mathcal{C}_t}(\tilde x_t) \right) + \nu_0^T(\tilde x_0-x_\mathrm{init}) \\
&& ~~~~~~ \nonumber + \sum_{t=0}^{N-1}\nu_{t+1}^T (F\tilde x_t-c_t) \\
&& ~~~~~~ \nonumber + \sum_{t=0}^{N-1}\Big( -\rho \tilde v_{t+1}^T (G_1\tilde x_t-\tilde z_{t+1}) + \frac{\rho}{2} \|G_1\tilde x_t-\tilde z_{t+1}\|_2^2 \Big) \\
&& ~~~~~~ \nonumber + \sum_{t=1}^N\Big( -\rho \tilde w_t^T (G_0\tilde x_t-\tilde z_t) + \frac{\rho}{2}\|G_0\tilde x_t-\tilde z_t\|_2^2	\Big).
\end{eqnarray}
The dual variables $\nu_t\in\reals^n, \; t=0,\ldots,N$ are associated with the equality constraints. We define the indicator function
for the polyhedral set $\mathcal C_t\subset \reals^{2n+m}$ as
\[
I_{\mathcal{C}_t}(\tilde x_t) =
\left\{ \begin{array}{ll} 0 & \tilde x_t \in \mathcal C_t\\
\infty & \mbox{otherwise}.
\end{array} \right.
\]
The augmented Lagrangian is minimized (and maximized) over the primal (and dual) variables in an iterative manner, for each
iteration of the algorithm. \\

By treating the dynamics and inequality constraints as explicit and minimizing with respect to $\tilde x_t, \; t=0,\ldots,N$,
we end up with the stage-wise QPs (\ref{eq:time_split_x0}), (\ref{eq:time_split_xt}), (\ref{eq:time_split_xT}). \\

The update of $\tilde z_t$ is given by solving
 \begin{eqnarray}
  \tilde z_t^{k+1} \nonumber &=& \argmin_{\tilde z_t} \Big( -\rho \tilde v_{t}^{kT}(G_1\tilde x_{t-1}^{k+1}-\tilde z_t) \\
               \nonumber &+& \frac{\rho}{2}\|G_1\tilde x_{t-1}^{k+1}-\tilde z_t\|_2^2 \\
               \nonumber &-& \rho \tilde w_t^{kT}(G_0\tilde x_t^{k+1}-\tilde z_t) +
  \frac{\rho}{2}\|G_0\tilde x_t^{k+1}-\tilde z_t\|_2^2 \Big) \\
               \nonumber &=& \argmin_{\tilde z_t} \Big( \frac{\rho}{2}\|\tilde z_t-G_1\tilde x_{t-1}^{k+1}+
                \tilde v_t^k\|_2^2 -
               \frac{\rho}{2}\|\tilde v_t^k\|_2^2 \\
               \nonumber &+& \frac{\rho}{2}\|\tilde z_t-G_0\tilde x_t^{k+1}+\tilde w_t^k\|_2^2 -
               \frac{\rho}{2}\|\tilde w_t^k\|_2^2 \Big) \Rightarrow \\
  \tilde z_t^{k+1} &=& \frac{G_0\tilde x_t^{k+1}+G_1\tilde x_{t-1}^{k+1}-\tilde v_t^k -\tilde w_t^k}{2},
                   \quad t=1,\ldots,N {\label{eq:time_split_z}}
 \end{eqnarray}

 Rearranging the terms in (\ref{eq:AugLag_time_split}) and completing the squares as before, the dual updates for
 $\tilde w_t$ and $\tilde v_t$ result by solving
 \begin{eqnarray*}
  \tilde w_t^{k+1} &=& \argmax_{\tilde w_t} \left( \frac{\rho}{2}\|\tilde w_t-G_0\tilde x_t^{k+1}+\tilde z_t^{k+1}\|_2^2 -
                                           \frac{\rho}{2}\|\tilde w_t\|_2^2 \right), \\
  \tilde v_t^{k+1} &=& \argmax_{\tilde v_t} \left( \frac{\rho}{2}\|\tilde v_t-G_1\tilde x_{t-1}^{k+1}+\tilde z_t^{k+1}\|_2^2 -
						\frac{\rho}{2}\|\tilde v_t\|_2^2 \right)
 \end{eqnarray*}
The dual updates are
\begin{eqnarray}
 \tilde w_t^{k+1} &=& \tilde w_t^k - G_0\tilde x_t^{k+1} + \tilde z_t^{k+1} {\label{eq:time_split_dual1}},\\
 \tilde v_t^{k+1} &=& \tilde v_t^k - G_1\tilde x_{t-1}^{k+1} + \tilde z_t^{k+1}, \quad t=1,\ldots,T.
 {\label{eq:time_split_dual2}}
\end{eqnarray}

Substituting (\ref{eq:time_split_dual1}) and (\ref{eq:time_split_dual2}) into (\ref{eq:time_split_z}) results in
$\tilde w_t^{k+1}+\tilde v_t^{k+1}=0$, for all $t$. Hence $\tilde z_t^{k+1}$ can be simplified to
\begin{equation}
 \tilde z_t^{k+1} \nonumber = \frac{G_0\tilde x_t^{k+1}+G_1\tilde x_{t-1}^{k+1}}{2}, \quad t=1,\ldots,N.
\end{equation}
This result shows that $\tilde z_t$ is the variable that enforces
consensus upon $x_t^{(t)}$ and $x_t^{(t-1)}$, the value of which is the average of the two.

\subsection*{Derivation of the three-set splitting algorithm updates}
The augmented Lagrangian for (\ref{eq:3sets_split_objective}) can be written as
\begin{eqnarray}
&& L_{\rho}(x,z,y,\tilde z,\tilde y,\nu) := \label{eq_AugLag_3sets_split} \\
&& ~~~~~~ \nonumber \frac{1}{2} x^{\I T}Mx^{\I} + q^Tx^{\I} + r + \nu^T(\tilde Ax^{\II}-\tilde b) \\
&& ~~~~~~ \nonumber - \rho \sum_{i=\I}^{\III}\tilde z^{iT} (x^{i}-z) + \frac{\rho}{2}\sum_{i=\I}^{\III}\|x^{i}-z\|_2^2 \\
&& ~~~~~~ \nonumber - \rho \tilde y^T(\tilde h - \tilde Hx^{\III}-y) + \frac{\rho}{2}\|\tilde h - \tilde Hx^{\III}-y\|_2^2 + I_{\{y \in \reals_+^p\}}.
\end{eqnarray}
As before, the dual variable $\nu$ is associated with the equality constraints. We define the indicator function
of the nonnegative orthant for a variable $x \in \reals$ as
\[
I_{\{x \in \reals_+\}} = x_+ =
\left\{ \begin{array}{ll} x & x \geq 0 \\
\infty & \mbox{otherwise}.
\end{array} \right.
\]
The function applies componentwise for vector variables.

The augmented Lagrangian is a smooth, quadratic function with respect to the variables $x^{\I}, x^{\II}, x^{\III}$, hence
the updates can be derived from taking the gradient equal to zero. This yields the solutions (\ref{eq:3_sets_x1}),
(\ref{eq:3_sets_x2}) and (\ref{eq:3_sets_x3}). \\

Grouping the terms that involve $z$ from (\ref{eq_AugLag_3sets_split}) we obtain
 \begin{eqnarray}
  z^{k+1} \nonumber &=& \argmin_{z} \left( \sum_{i=\I}^{\III} \frac{\rho}{2}\|z - \left(x^{i}\right)^{k+1} +
                        \left(\tilde z^{i}\right)^k\|^2 \right) \Rightarrow \\
  z^{k+1} &=& \frac{1}{3}\sum_{i=\I}^{\III} \left(\left(x^{i}\right)^{k+1} - \left(\tilde z^{i}\right)^k\right).
                     {\label{eq:3sets_split_z}}
 \end{eqnarray}
 Maximizing (\ref{eq_AugLag_3sets_split}) over $\tilde z^{i}, i=\I, \II, \III$ results in
 \[
  \left(\tilde z^{i}\right)^{k+1} = \argmax_{\tilde z^{i}} \left( -\frac{\rho}{2}\|\tilde z^{i}\|_2^2 +
                                                         \frac{\rho}{2}
                                                         \|\tilde z^{i}-\left(x^{i}\right)^{k+1}+z^{k+1}\|_2^2 \right).
 \]

Hence, the update for $\tilde z^{i}, i=\I, \II, \III$ is
\begin{equation}\label{eq:3sets_split_dual2}
 \left(\tilde z^{i}\right)^{k+1} = \left(\tilde z^{i}\right)^k - \left(x^{i}\right)^{k+1} + z^{k+1}, \; i=\I, \II, \III.
\end{equation}

Similarly, the update for $\tilde y$ is
\begin{equation*}
 \tilde y^{k+1} = \tilde y^k + y^{k+1} - h + H\left(x^{\III}\right)^{k+1}.
\end{equation*}

Taking the sum of the dual variables $\tilde z^{i}, \; i=\I, \II, \III$ from (\ref{eq:3sets_split_dual2}) and substituting into
(\ref{eq:3sets_split_z}) simplifies the $z$-update and we obtain
\[
 z^{k+1} = \frac{1}{3}\sum_{i=\I}^{\III} \left(x^{i}\right)^{k+1}.
\]
Hence, the global variable $z$ is the average of the three estimates of $x$.

\section*{ACKNOWLEDGMENTS}
The main part of this work was carried out at Stanford. The authors would like to thank Stephen Boyd and Brendan O'Donoghue for helpful discussions.


\bibliographystyle{ieeetr}
\bibliography{ECC13refs}

\end{document}